\input amstex
\input amsppt.sty
\magnification=\magstep1
\vsize=22.2truecm
\baselineskip=16truept
\NoBlackBoxes
\nologo
\pageno=1
\topmatter

\def\Proof{\noindent{\it Proof}}

\def\Ack{\medskip\noindent {\bf Acknowledgment}}
\def\pmod #1{\ (\roman{mod}\ #1)}
\title  On a generalization of Carlitz's congruence\endtitle
\author
Hao Pan
\endauthor
\address
Department of Mathematics, Nanjing University,
Nanjing 210093, People's Republic of China
\endaddress
\email{haopan79\@yahoo.com.cn}\endemail
\abstract Let $p$ be an odd prime and $a$ be a positive integer. We
show that
$$
\sum_{k=0}^{p-1}(-1)^{(a-1)k}\binom{p-1}{k}^a\equiv2^{a(p-1)}+\frac{a(a-1)(3a-4)}{48}p^3B_{p-3}\pmod{p^4},
$$
which is a generalization of a congruence due to Carlitz.
\endabstract
\endtopmatter
\document
\TagsOnRight
\heading
1. Introduction
\endheading
As early as 1895, with the help of De Moivre's theorem, Morley
[13] (or cf. [9]) proved a beautiful congruence for binomial
coefficients:
$$
(-1)^{(p-1)/2}\binom{p-1}{(p-1)/2}\equiv 4^{p-1}\pmod{p^3}\tag 1.1
$$
for any prime $p\geqslant 5$. And Carlitz [3] extended Morley's congruence as follows:
$$
(-1)^{(p-1)/2}\binom{p-1}{(p-1)/2}\equiv 4^{p-1}+\frac{1}{12}p^3B_{p-3}\pmod{p^4}\tag 1.2
$$
for each odd prime $p$, where $B_n$ are the Bernoulli numbers given by
$$
\frac{t}{e^t-1}=\sum_{n=0}^{\infty}\frac{B_n}{n!}t^n.
$$

On the other hand, some combinatorial and arithmetical properties of the binomial sums
$$
\sum_{k=0}^n\binom{n}{k}^a\text{ and }\sum_{k=0}^n(-1)^k\binom{n}{k}^a
$$
have been investigated by several authors (e.g., Calkin [4], Cusick [5], McIntosh [11], Perlstadt [13]). Indeed, we know [8, Eqs. (3.81) and (6.6)]that
$$
\sum_{k=0}^n(-1)^k\binom{2n}{k}^2=(-1)^n\sum_{k=0}^n\binom{n}{k}^2=(-1)^n\binom{2n}{n}\tag 1.3
$$
and
$$
\sum_{k=0}^n(-1)^k\binom{2n}{k}^3=(-1)^n\binom{2n}{n}\binom{3n}{n}.\tag 1.4
$$
However, by using asymptotic methods, de Bruijn [1] has showed that no closed form exists for the sum $\sum_{k=0}^n(-1)^k\binom{n}{k}^a$ when $a\geqslant 4$. And
Wilf proved (in a personal communication with Calkin, see [4]) that the sum $\sum_{k=0}^n\binom{n}{k}^a$ has no closed form provided that $3\leqslant a\leqslant 9$.

Recently Chamberland and Dilcher [6] studied the congruences for the sum
$$
u_{a,b}^{\epsilon}(n)=\sum_{k=0}^{2n}(-1)^{\epsilon k}\binom{n}{k}^a\binom{2n}{k}^b
$$
where $a, b\geqslant 0$ and $\epsilon\in\{0,1\}$. For example, they proved that for any prime $p\geqslant5$
$$
u_{a,b}^{\epsilon}(p)\equiv1+(-1)^{\epsilon}2^b\pmod{p^3},
$$
unless $(\epsilon,a,b)=(0,0,1)$ or $(0,1,0)$. In fact, the sum $\sum_{k=0}^{p-1}(-1)^{(a-1)k}\binom{p-1}{k}^a$ has been considered by Cai and Granville in [2, Theorem 6]. And
they proved that
$$
\sum_{k=0}^{p-1}(-1)^{(a-1)k}\binom{p-1}{k}^a\equiv2^{a(p-1)}\pmod{p^3}
$$
for a prime $p\geqslant 5$. In view of (1.3), Cai and Granville's result generalized the congruence of Morley.
In this paper we shall give such a generalization of of Carlitz's congruence (1.2).
\proclaim{Theorem 1.1} Let $p$ be an odd prime and $a$ be a positive integer. Then
$$
\sum_{k=0}^{p-1}(-1)^{(a-1)k}\binom{p-1}{k}^a\equiv2^{a(p-1)}+\frac{a(a-1)(3a-4)}{48}p^3B_{p-3}\pmod{p^4}.\tag 1.5
$$
\endproclaim
There are another simple consequences of Theorem 1.1.
\proclaim{Corollary 1.2} Let $p$ be an odd prime.
Then
$$
\sum_{k=0}^{p-1}\binom{p-1}{k}^3\equiv8^{p-1}+\frac{5}{8}p^3B_{p-3}\pmod{p^4},\tag 1.6
$$
$$
\sum_{k=0}^{p-1}(-1)^k\binom{p-1}{k}^4\equiv16^{p-1}+2p^3B_{p-3}\pmod{p^4}\tag 1.7
$$
and
$$
\sum_{k=0}^{p-1}\binom{p-1}{k}^5\equiv32^{p-1}+\frac{55}{12}p^3B_{p-3}\pmod{p^4}.\tag 1.8
$$
\endproclaim
\heading
2. Proof of Theorem 1.1
\endheading
The Bernoulli polynomials $B_n(x)$ are defined by
$$
\frac{te^{xt}}{e^t-1}=\sum_{n=0}^{\infty}\frac{B_n(x)}{n!}t^n.
$$
Clearly $B_n=B_n(0)$. Also we have
$$
\sum_{k=1}^{n-1}k^{m-1}=\frac{B_m(n)-B_m}{m}
$$
for any positive integers $n$ and $m$. For more properties of Bernoulli numbers and Bernoulli polynomials, the readers may refer to [7] and [10].

\proclaim{Lemma 2.1} Let $p\geqslant 5$ be a prime. Then

\noindent{\rm (i)}
$$
\sum_{k=1}^{(p-1)/2}\frac{1}{k}\equiv-2q_2(p)+pq_2^2(p)-\frac{2}{3}p^2q_2^3(p)-\frac{7}{12}p^2B_{p-3}\pmod{p^3},\tag 2.1
$$
where $q_2(p)=(2^{p-1}-1)/p$.

\noindent{\rm (ii)}
$$
\sum_{k=1}^{(p-1)/2}\frac{1}{k^n}\equiv\cases
\frac{n(2^{n+1}-1)}{2(n+1)}pB_{p-n-1}\pmod{p^2}&\text{if }2\mid n,\\
-\frac{2(2^{n-1}-1)}{n}B_{p-n}\pmod{p}&\text{if }2\nmid n
\endcases\tag 2.2
$$
for $2\leqslant n\leqslant p-2$.
\endproclaim
\Proof. See Theorem 5.2 and Corollary 5.2 in Sun's paper [14].\qed
\proclaim{Lemma 2.2} Let $p\geqslant 5$ be a prime. Then we have
$$
\sum_{\Sb 1\leqslant j<k\leqslant p-1\\2\mid k\endSb}\frac{1}{jk^2}\equiv\frac{5}{8}B_{p-3}\pmod{p}\tag 2.4
$$
and
$$
\sum_{\Sb 1\leqslant j<k\leqslant p-1\\2\mid k\endSb}\frac{1}{j^2k}\equiv-\frac{3}{8}B_{p-3}\pmod{p}.\tag 2.5
$$
\endproclaim
\Proof. It follows from Lemma 2.1 that $\sum_{k=1}^{(p-1)/2}k^n$ is divisible by $p$ if $2\leqslant n\leqslant p-3$ and $n$ is even.
And we know that $B_n=0$ for all odd $n\geqslant 3$. Hence
$$
\align
\sum_{\Sb 1\leqslant j<k\leqslant p-1\\2\mid k\endSb}\frac{1}{jk^2}\equiv&\sum_{\Sb 1\leqslant j<k\leqslant p-1\\2\mid k\endSb}\frac{j^{p-1}}{k^2}=\sum_{\Sb 1\leqslant k\leqslant p-1\\2\mid k\endSb}\frac{B_{p-1}(k)-B_{p-1}}{(p-1)k^2}\\
=&\sum_{\Sb 1\leqslant k\leqslant p-1\\2\mid k\endSb}\sum_{i=1}^{p-1}\binom{p-1}{i}\frac{k^{i-2}B_{p-1-i}}{p-1}\equiv\binom{p-1}{2}\frac{B_{p-3}}{2}-\sum_{\Sb k=1\\2\mid k\endSb}^{p-1}\frac{k^{p-4}}{2}\\
\equiv&\binom{p-1}{2}\frac{B_{p-3}}{2}-\sum_{\Sb k=1\\2\mid k\endSb}^{p-1}\frac{1}{2k^3}\equiv\frac{5}{8}B_{p-3}\pmod{p}.
\endalign
$$
This concludes the proof of (2.4). And we left the proof of (2.5) as an exercise for the readers.
\qed

\medskip\noindent{\it Proof of Theorem 1.1}.
Assume that $p\geqslant 5$. For any $1\leqslant r<p$, we have
$$
(-1)^r\binom{p-1}{r}=\prod_{k=1}^r\frac{k-p}{k}\equiv1-\sum_{k=1}^r\frac{p}{k}+\sum_{1\leqslant j<k\leqslant r}\frac{p^2}{jk}-\sum_{1\leqslant i<j<k\leqslant r}\frac{p^3}{ijk}\pmod{p^4}.
$$
Therefore
$$
\align
(-1)^{ar}\binom{p-1}{r}^a
\equiv&1-a\big(\sum_{k=1}^r\frac{p}{k}-\sum_{1\leqslant j<k\leqslant r}\frac{p^2}{jk}+\sum_{1\leqslant i<j<k\leqslant r}\frac{p^3}{ijk}\big)+\binom{a}{2}\big(\sum_{k=1}^r\frac{p}{k}\big)^2\\
&-2\binom{a}{2}\big(\sum_{k=1}^r\frac{p}{k}\big)\big(\sum_{1\leqslant j<k\leqslant r}\frac{p^2}{jk}\big)-\binom{a}{3}\big(\sum_{k=1}^r\frac{p}{k}\big)^3\pmod{p^4}.
\endalign
$$
Note that
$$
\big(\sum_{k=1}^r\frac{p}{k}\big)^2=2\sum_{1\leqslant j<k\leqslant r}\frac{p^2}{jk}+\sum_{k=1}^r\frac{p^2}{k^2}.
$$
Also it is easy to check that
$$
\big(\sum_{i=1}^r\frac{p}{i}\big)\big(\sum_{1\leqslant j<k\leqslant r}\frac{p^2}{jk}\big)=3\sum_{1\leqslant i<j<k<i\leqslant r}\frac{p^3}{ijk}+\sum_{1\leqslant j<k\leqslant r}\frac{p^3}{j^2k}+\sum_{1\leqslant j<k\leqslant r}\frac{p^3}{jk^2}
$$
and
$$
\big(\sum_{k=1}^r\frac{p}{k}\big)^3=6\sum_{1\leqslant i<j<k<i\leqslant r}\frac{p^3}{ijk}+3\sum_{1\leqslant j<k\leqslant r}\frac{p^3}{j^2k}+3\sum_{1\leqslant j<k\leqslant r}\frac{p^3}{jk^2}+\sum_{k=1}^r\frac{p^3}{k^3}.
$$
Thus
$$
\align
&\sum_{r=1}^{p-1}(-1)^{(a-1)r}\binom{p-1}{r}^a\\
\equiv&\sum_{r=1}^{p-1}(-1)^{r}\bigg(1-\sum_{k=1}^r\frac{ap}{k}+\sum_{1\leqslant j<k\leqslant r}\frac{a^2p^2}{jk}-\sum_{1\leqslant i<j<k\leqslant r}\frac{a^3p^3}{ijk}
+\binom{a}{2}\sum_{k=1}^r\frac{p^2}{k^2}\\
&-\binom{a}{2}\big(\sum_{1\leqslant j<k\leqslant r}\frac{ap^3}{j^2k}+\sum_{1\leqslant j<k\leqslant r}\frac{ap^3}{jk^2}\big)-\binom{a}{3}\sum_{k=1}^r\frac{p^3}{k^3}\bigg)\\
=&\bigg(-\sum_{k=1}^{p-1}\frac{ap}{k}+\sum_{1\leqslant j<k\leqslant p-1}\frac{a^2p^2}{jk}-\sum_{1\leqslant i<j<k\leqslant p-1}\frac{a^3p^3}{ijk}
+\binom{a}{2}\sum_{k=1}^{p-1}\frac{p^2}{k^2}\\
&-\binom{a}{2}\big(\sum_{1\leqslant j<k\leqslant p-1}\frac{ap^3}{j^2k}+\sum_{1\leqslant j<k\leqslant p-1}\frac{ap^3}{jk^2}\big)-\binom{a}{3}\sum_{k=1}^{p-1}\frac{p^3}{k^3}\bigg)\sum_{r=k}^{p-1}(-1)^{r}\\
=&-\sum_{\Sb k=1\\2\mid k\endSb}^{p-1}\frac{ap}{k}+\sum_{\Sb 1\leqslant j<k\leqslant p-1\\2\mid k\endSb}\frac{a^2p^2}{jk}-\sum_{\Sb 1\leqslant i<j<k\leqslant p-1\\2\mid k\endSb}\frac{a^3p^3}{ijk}
+\binom{a}{2}\sum_{\Sb k=1\\2\mid k\endSb}^{p-1}\frac{p^2}{k^2}\\
&-\binom{a}{2}\big(\sum_{\Sb 1\leqslant j<k\leqslant p-1\\2\mid k\endSb}\frac{ap^3}{j^2k}+\sum_{\Sb 1\leqslant j<k\leqslant p-1\\2\mid k\endSb}\frac{ap^3}{jk^2}\big)-\binom{a}{3}\sum_{\Sb k=1\\2\mid k\endSb}^{p-1}\frac{p^3}{k^3}\pmod{p^4}.
\tag 2.6
\endalign
$$
Now we only need to determine
$$
\sum_{\Sb 1\leqslant j<k\leqslant p-1\\2\mid k\endSb}\frac{1}{jk}\pmod{p^2}\text{ and }\sum_{\Sb 1\leqslant i<j<k\leqslant p-1\\2\mid k\endSb}\frac{1}{ijk}\pmod{p}.
$$
Letting $a=1$ in (2.6), we obtain that
$$
2^{p-1}-1\equiv-\sum_{\Sb k=1\\2\mid k\endSb}^{p-1}\frac{p}{k}+\sum_{\Sb 1\leqslant j<k\leqslant p-1\\2\mid k\endSb}\frac{p^2}{jk}-\sum_{\Sb 1\leqslant i<j<k\leqslant p-1\\2\mid k\endSb}\frac{p^3}{ijk}\pmod{p^4},
$$
whence
$$
\sum_{\Sb 1\leqslant j<k\leqslant p-1\\2\mid k\endSb}\frac{1}{jk}\equiv\sum_{\Sb 1\leqslant i<j<k\leqslant p-1\\2\mid k\endSb}\frac{p}{ijk}+\frac{1}{2}q_2^2(p)-\frac{1}{3}pq_2^3(p)-\frac{7}{24}pB_{p-3}\pmod{p^2}.\tag 2.7
$$
Also setting $a=2$ in (2.6), then by Carlitz's congruence (1.2),
$$
\align
&\sum_{\Sb 1\leqslant j<k\leqslant p-1\\2\mid k\endSb}\frac{4p^2}{jk}-\sum_{\Sb k=1\\2\mid k\endSb}^{p-1}\frac{2p}{k}-\sum_{\Sb 1\leqslant i<j<k\leqslant p-1\\2\mid k\endSb}\frac{8p^3}{ijk}
+\sum_{\Sb k=1\\2\mid k\endSb}^{p-1}\frac{p^2}{k^2}-\sum_{\Sb 1\leqslant j<k\leqslant p-1\\2\mid k\endSb}\big(\frac{2p^3}{j^2k}+\frac{2p^3}{jk^2}\big)\\
\equiv&\sum_{\Sb 1\leqslant j<k\leqslant p-1\\2\mid k\endSb}\frac{4p^2}{jk}-\sum_{\Sb 1\leqslant i<j<k\leqslant p-1\\2\mid k\endSb}\frac{8p^3}{ijk}+2pq_2(p)-p^2q_2^2(p)+\frac{2}{3}p^3q_2^3(p)+\frac{2}{3}p^3B_{p-3}\\
\equiv&(-1)^{(p-1)/2}\binom{p-1}{(p-1)/2}-1\equiv4^{p-1}+\frac{1}{12}p^3B_{p-3}-1\pmod{p^4},
\endalign
$$
that is,
$$
\sum_{\Sb 1\leqslant j<k\leqslant p-1\\2\mid k\endSb}\frac{1}{jk}\equiv\sum_{\Sb 1\leqslant i<j<k\leqslant p-1\\2\mid k\endSb}\frac{2p}{ijk}+\frac{1}{2}q_2^2(p)-\frac{1}{6}pq_2^3(p)-\frac{7}{48}pB_{p-3}\pmod{p^2}.\tag 2.8
$$
Combining (2.7) and (2.8), we have
$$
\sum_{\Sb 1\leqslant j<k\leqslant p-1\\2\mid k\endSb}\frac{1}{jk}\equiv\frac{1}{2}(q_2^2(p)-pq_2^3(p))-\frac{7}{16}pB_{p-3}\pmod{p^2}\tag 2.9
$$
and
$$
\sum_{\Sb 1\leqslant j<k\leqslant p-1\\2\mid k\endSb}\frac{1}{ijk}\equiv-\frac{1}{6}q_2^3(p)-\frac{7}{48}B_{p-3}\pmod{p}.\tag 2.10
$$
Substituting (2.9) and (2.10) in (2.6), it follows that
$$
\align
&\sum_{r=0}^{p-1}(-1)^{(a-1)r}\binom{p-1}{r}^a\\
\equiv&1+\binom{a}{1}pq_2(p)+\binom{a}{2}p^2q_2^2(p)+\binom{a}{3}p^3q_2^3(p)+\big(\frac{1}{12}\binom{a}{2}+\frac{3}{8}\binom{a}{3}\big)p^3B_{p-3}\\
\equiv&\sum_{j=0}^{a}\binom{a}{j}(2^{p-1}-1)^j+\frac{a(a-1)(3a-4)}{48}p^3B_{p-3}\pmod{p^4}.
\endalign
$$

Finally, when $p=3$,
$$
\align
2^{2a}-\sum_{k=0}^{2}(-1)^{(a-1)k}\binom{2}{k}^a=&(3+1)^{a}-(2-(-1)^a(3-1)^a)\\
\equiv&\sum_{j=0}^{3}\binom{a}{j}3^j+\sum_{j=0}^{3}\binom{a}{j}(-3)^j-2=9a(a-1)\pmod{3^4}.
\endalign
$$
And
$$
\frac{a(a-1)(3a-4)}{48}\cdot3^3\cdot B_{3-3}+9a(a-1)=\frac{27}{16}a(a-1)(a+4)\equiv0\pmod{3^4}.
$$
All are done.

\Ack. I thank Professor Tian-Xin Cai for informing me the result of (1.5) modulo $p^3$.
I also thank my advisor, Professor Zhi-Wei Sun, for his help on this paper.

\enddocument